\documentclass[12pt,leqno]{article}
\usepackage{amsthm,amsfonts,amssymb,amsmath,eufrak,oldgerm}
\usepackage{epsfig}
\numberwithin{equation}{section} 

\newcommand{\des}{{\alpha_*}}
\renewcommand\d{\partial}

\renewcommand\b{\beta}

\renewcommand\o{\omega}
\newcommand\R{\mathbb R}
\newcommand\C{\mathbb C}

\def\eps{\varepsilon}
\def\e{\varepsilon}

\newcommand\br{\begin{remark}}
\newcommand\er{\end{remark}}
\newcommand\bp{\begin{pmatrix}}
\newcommand\ep{\end{pmatrix}}
\newcommand\be{\begin{equation}}
\newcommand\ee{\end{equation}}
\newcommand\ba{\begin{equation}\begin{aligned}}
\newcommand\ea{\end{aligned}\end{equation}}

\newcommand{\bap}{\begin{app}}
\newcommand{\eap}{\end{app}}
\newcommand{\begs}{\begin{exams}}
\newcommand{\eegs}{\end{exams}}
\newcommand{\beg}{\begin{example}}
\newcommand{\eeg}{\end{exaplem}}
\newcommand{\bpr}{\begin{proposition}}
\newcommand{\epr}{\end{proposition}}
\newcommand{\bt}{\begin{theorem}}
\newcommand{\et}{\end{theorem}}
\newcommand{\bc}{\begin{corollary}}
\newcommand{\ec}{\end{corollary}}
\newcommand{\bl}{\begin{lemma}}
\newcommand{\el}{\end{lemma}}
\newcommand{\bd}{\begin{definition}}
\newcommand{\ed}{\end{definition}}
\newcommand{\brs}{\begin{remarks}}
\newcommand{\ers}{\end{remarks}}

\newcommand{\Id}{{\rm Id }}

\newcommand{\Res}{{\rm Residue}}

\newtheorem{theorem}{Theorem}[section]
\newtheorem{proposition}[theorem]{Proposition}
\newtheorem{corollary}[theorem]{Corollary}
\newtheorem{lemma}[theorem]{Lemma}
\newtheorem{definition}[theorem]{Definition}

\newtheorem{example}[theorem]{Example}
\newtheorem{remark}[theorem]{Remark}

\newtheorem{exams}[theorem]{Examples}

\newcommand\cB{{\cal  B}}

\newcommand\cW{{\cal  W}}
\newcommand\cZ{{\cal  Z}}

\newcommand\cG{{\cal  G}}

\newcommand\cN{{\cal  N}}
\newcommand\cE{{\cal  E}}
\newcommand\cF{{\cal  F}}

\newcommand\cM{{\mathcal M}}
\newcommand\cT{{\mathcal T}}

\title{
Center stable manifolds for quasilinear parabolic pde and
conditional stability of nonclassical viscous shock waves
}

\author{\sc \small 
Kevin Zumbrun\thanks{Indiana University, Bloomington, IN 47405;
kzumbrun@indiana.edu:
Research of K.Z. was partially supported
under NSF grants no. DMS-0300487 and DMS-0801745.
 }}
\begin{document}

\maketitle


\begin{abstract}
Motivated by the study of conditional stability of traveling
waves, we give an elementary $H^2$ center stable manifold construction for
quasilinear parabolic PDE, 
sidestepping apparently delicate regularity issues by the 
combination of a carefully chosen implicit fixed-point scheme
and straightforward
time-weighted $H^s$ energy estimates. 
As an application, we show conditional stability of Lax- or
undercompressive shock waves of general
quasilinear parabolic systems of conservation
laws by a pointwise stability analysis on the center stable manifold.
\end{abstract}

\tableofcontents
\bigbreak
\section{Introduction}
In this paper, extending our previous work in the semilinear case \cite{Z5},
we show by an elementary argument
that an asymptotically constant stationary solution
\be\label{soln}
u(x,t)\equiv \bar u(x),
\quad
|\bar u(x)-u_\pm|\le Ce^{-\theta|x|},
\ee
$\theta>0$, of a general quasilinear second-order parabolic system
\be\label{rcd}
u_t=\cF(u):= b(u)u_{xx} - h(u,u_x),
\ee
$x$, $t\in \R$, $u$, $h\in \R^n$, $b\in \R^{n\times n}$,
possesses a local center stable manifold with respect to $H^2$.
Combining this result with ideas from \cite{ZH,HZ,RZ,TZ2,TZ3,Z5}
we then establish conditional stability of nonclassical viscous shock 
solutions of strictly parabolic systems of conservation laws, similarly
as was done in \cite{Z5} for Lax shocks of semilinear parabolic systems.

\subsection{Existence of center stable manifolds}
We briefly describe our results.
Following \cite{Z5}, assume:

\medbreak
(A0) $h \in C^{k+1}$, $k\ge 2$.

(A1) $\Re \sigma(b)\ge \theta>0$.

(A2) The linearized operator $L=\frac{\d \cF}{\d u}(\bar u)$ about
$\bar u$ has $p$ unstable (positive real part)
eigenvalues, with the rest of its spectrum  of nonpositive real part.

(A3) $|\partial_x^j \bar u (x)|\le Ce^{-\theta |x|}$, $\theta>0$,
for $1\le j\le k+2$.

\medbreak
Then, we have the following basic version of
the Center Stable Manifold Theorem.

\bpr \label{t:cspde}
Under assumptions (A0)--(A2), there exists 
in an $H^2$ neighborhood of $\bar u$
a Lipschitz (with respect to $H^2$) 
center stable manifold $\cM_{cs}$,
tangent to quadratic order at $\bar u$ to the center stable subspace 
$\Sigma_{cs}$ of $L$ in the sense that
\be\label{tanbd1}
|\Pi_{u}(u-\bar u)|_{H^2}\le C |\Pi_{cs}(u-\bar u)|_{H^2}^2
\ee
for $u\in \cM_{cs}$ where $\Pi_{cs}$ and $\Pi_u$ denote the center-stable
and stable eigenprojections of $L$,
that is (locally) invariant under the forward time-evolution of \eqref{rcd}
and contains all solutions that remain sufficiently
close to $\bar u$ in forward time.  In general it is not unique.
\epr

Combining with ideas of \cite{TZ1,Z5}, we readily obtain by the 
same technique the following
improved version respecting the underlying translation-invariance
of \eqref{rcd}, a property that is important in the applications \cite{Z5}.

\bt \label{t:cstrans}
Under assumptions (A0)--(A2), there exists 
in an $H^2$ neighborhood of the set of translates of $\bar u$
a translation invariant Lipschitz (with respect to $H^2$) 
center stable manifold $\cM_{cs}$,
tangent to quadratic order at $\bar u$ to the center stable subspace 
$\Sigma_{cs}$ of $L$ in the sense of \eqref{tanbd1},
that is (locally) invariant under the forward time-evolution of \eqref{rcd}
and contains all solutions that remain bounded and sufficiently
close to a translate of $\bar u$ in forward time.  In general it is not unique.
\et

\subsection{Conditional stability of nonclassical shocks}
Now, specialize to the case
\be\label{shockcase}
h(u,u_x)=f(u)_x - (db(u)u_x)u_x
\ee
that \eqref{rcd} corresponds to
a parabolic system of conservation laws in standard form
\be\label{cons}
u_t+ f(u)_{x}= (b(u)u_{x})_x,
\ee
$u,\, f\in \R^n$, $b\in \R^{n\times n}$, 
$x\,,t\in \R$.  Following \cite{ZH, HZ}, assume: 

(H0) $f$, $b \in C^{k+2}$, $k\ge 2$.

(H1) $\Re \sigma(b)\ge \theta>0$.

(H2) $A_\pm:=df(u_\pm)$ have simple, real, nonzero eigenvalues.

(H3) $\Re \sigma(df i\xi-b|\xi|^2)(u_\pm)\le -\theta |\xi|^2$, $\theta>0$,
for all $\xi\in \R$.

(H4) Nearby $\bar u$, the set of all solutions
\eqref{soln} connecting the same values $u_\pm$
forms a smooth manifold $\{\bar u^\alpha\}$, $\alpha\in \mathcal{U}\subset 
\mathbb{R}^\ell$, $\bar u^0=\bar u$.

(H5) The dimensions of the unstable subspace of $df(u_-)$
and the stable subspace of $df(u_+)$ sum to either $n+\ell$, $\ell=1$
(pure Lax case), $n+\ell$, $\ell>1$ (pure overcompressive case),
or $\le n$ with $\ell=1$ (pure undercompressive case).
\medbreak

Assume further the following {\it spectral genericity} conditions.
\medbreak

(D1) $L$ has no nonzero imaginary eigenvalues. 

(D2) The Evans function $D(\lambda)$ associated with $L$
vanishes at $\lambda=0$ to precisely order $\ell$.

\medbreak
Here, the {\it Evans function}, as defined in \cite{ZH,GZ}
denotes a certain Wronskian associated with eigenvalue ODE $(L-\lambda)w=0$,
whose zeros correspond in location and multiplicity with the eigenvalues
of $L$.  For history and basic properties of the Evans function, 
see, e.g. \cite{AGJ,PW,GZ,MaZ1} and references therein.

As discussed in \cite{ZH,MaZ1}, (D2) corresponds in the
absence of a spectral gap to a generalized notion of simplicity of the 
embedded eigenvalue $\lambda=0$ of $L$.
Thus, (D1)--(D2) together correspond to the assumption that there are
no additional (usual or generalized) eigenvalues on the imaginary
axis other than the transational eigenvalue at $\lambda=0$;
that is, the shock is not in transition between different degrees of
stability, but has stability properties that are insensitive to 
small variations in parameters.

With these assumptions, we obtain our remaining results
characterizing the stability properties of $\bar u$.
Denoting by
\begin{equation}\label{aj}
a_1^\pm<a_2^\pm < \cdots < a_n^\pm
\end{equation}
the eigenvalues of the limiting convection matrices $A_\pm:= df(u_\pm)$,
define
\begin{equation}\label{theta}
\theta(x,t):=
\sum_{a_j^-<0}(1+t)^{-1/2}e^{-|x-a_j^-t|^2/Mt}
+ \sum_{a_j^+>0}(1+t)^{-1/2}e^{-|x-a_j^+t|^2/Mt},
\end{equation}
\begin{equation}\label{psi1}
\begin{aligned}
\psi_1(x,t)&:=
\chi(x,t)\sum_{a_j^-<0}
(1+|x|+t)^{-1/2} (1+|x-a_j^-t|)^{-1/2}\\
&\quad+
\chi(x,t)\sum_{a_j^+>0}
(1+|x|+t)^{-1/2} (1+|x-a_j^+t|)^{-1/2},\\
\end{aligned}
\end{equation}
and
\begin{equation}\label{psi2}
\begin{aligned}
\psi_2(x,t)&:=
(1-\chi(x,t)) (1+|x-a_1^-t|+t^{1/2})^{-3/2}\\
&\quad +(1-\chi(x,t)) (1+|x-a_n^+t|+t^{1/2})^{-3/2},
\end{aligned}
\end{equation}
where $\chi(x,t)=1$ for $x\in [a_1^-t, a_n^+t]$ and zero otherwise,
and $M>0$ is a sufficiently large constant.

\bpr\label{hA}
Conditions (H0)--(H5) imply (A0)--(A3), so that there
exists a translation-invariant center stable manifold $\cM_{cs}$
of $\bar u$ and its translates.
\epr

\bt\label{t:mainstab}
Under (H0)--(H5) and (D1)--(D2), $\bar u$ is nonlinearly
phase-asymptotically orbitally stable under sufficiently 
small perturbations $v_0\in H^4$ lying on the codimension $p$ center stable 
manifold $\cM_{cs}$ of $\bar u$ and its translates
with $|(1+|x|^2)^{3/4}v_0(x)|_{H^4}\le E_0$ sufficiently small,
where $p$ is the number of unstable eigenvalues of $L$,
in the sense that, for some $\alpha(\cdot)$, $\des$,
\begin{equation} \label{pointwise}
\begin{aligned}
|\partial_x^r\big(u(x,t)-\bar u^{\des + \alpha(t)}(x)\big)|&\le C E_0
(\theta+\psi_1+\psi_2)(x,t),\\
|u(\cdot,t)-\bar u^{\des + \alpha(t)}|_{H^4}&\le C E_0(1+t)^{-\frac{1}{4}},\\
|\des|&\le CE_0,\\
  |\alpha(t)|&\le C E_0(1+t)^{-1/2},\\
  |\dot \alpha (t)|&\le C E_0 (1+t)^{-1},\\
\end{aligned}
\end{equation}
$0\le r\le 1$,
where $u$ denotes the solution of \eqref{cons} with initial
data $u_0=\bar u+v_0$ and $\bar u^\alpha$ is as in (H4).
Moreover, $\bar u$ is orbitally unstable with respect to small $H^2$ 
perturbations not lying in $\cM$, in the sense that the 
corresponding solution leaves a fixed-radius neighborhood of the 
set of translates of $\bar u$ in finite time.
\et

\begin{remark}\label{rate}
\textup{
Pointwise bound \eqref{pointwise} yields as a corollary
the sharp $L^p$ decay rate
\begin{equation}
\label{Lp}
|u(x,t)-\bar u^{\des+\alpha(t)}(x)|_{L^p}\le C E_0
(1+t)^{-\frac{1}{2}(1-\frac{1}{p})}, \quad 1\le p\le \infty
\end{equation}
obtained by the $L^p$ (rather than pointwise) analysis of \cite{Z5}.
However, we obtain here the additional information that phase $\des+\alpha$
approaches a limit time-asymptotically at rate $t^{-\frac{1}{2}}$.
}
\end{remark}

\br
\textup{
Mixed over- or undercompressive type shocks are also possible
\cite{ZH,HZ}, though we do not know of any physical examples.
These can be treated with further effort as described in \cite{HZ}, 
Sections 5--6, or \cite{RZ}.
}
\er

\subsection{Discussion and open problems}\label{discuss}
As discussed in \cite{GJLS}, for semilinear problems like those
considered in \cite{GJLS,Z5},
for which the nonlinear part of the associated evolution equations consists
of a relatively compact perturbation of the linear part,
construction of invariant manifolds reduces essentially to verification
of a {spectral mapping theorem} for the linearized flow, after which
the construction follows by already-well-developed theory.
However, for quasilinear equations, 
the usual fixed-point construction of the standard theory
does not close due to apparent {loss of regularity}.
This appears to be a general difficulty in the construction of
invariant manifolds for quasilinear systems;
see \cite{Li,LPS1,LPS2} for further discussion.

We overcome this in the present case in a very simple way, by
(i) introducing a carefully chosen implicit fixed-point scheme for which the
infinite-dimensional center stable part satisfies a standard
Cauchy problem forced by the finite-dimensional (hence harmless
in terms of regularity) unstable part, and
(ii) making use of a straightforward 
time-weighted 
energy estimate
for the Cauchy problem to control higher derivatives by lower ones,
the latter of which may be estimated in standard linear fashion.
For related arguments, see for example \cite{MaZ2,MaZ3,Z2,Z3,RZ,TZ3}.

Existence of invariant manifolds of
quasilinear parabolic systems has been treated by
quite different methods in \cite{LPS1,LPS2} via
a detailed study of the smoothing properties of the linearized flow.
The advantage of the present method, besides its simplicity, is that
carries over in straightforward fashion to the 
compressible Navier--Stokes 
equations of gas dynamics and MHD \cite{Z6} for which the linearized
flow, being of hyperbolic--parabolic type,
is only partially smoothing.
A disadvantage is that we do not see by this technique
how to obtain smoothness of the center stable manifold,
but only Lipschitz continuity; this seems to be the
price of our simple energy-based approach.

As we are mainly interested in stability, it is not smoothness
but quadratic-order tangency \eqref{tanbd1} that is our main concern.
An observation of possibly general use is that
this weaker property is satisfied whenever the underlying flow is 
$C^2$ at $\bar u$, whereas global $C^2$ regularity 
of the center stable manifold would require
global $C^{2+\alpha}$ regularity, $\alpha>0$, of the flow. 
Moreover, it is easily verified in the course of the 
standard construction of a Lipschitz invariant manifold, and so
we obtain this key property with essentially no extra effort.

Regarding conditional stability, the main novelty 
of the present analysis is that we carry out a pointwise iteration 
scheme in order to
treat nonclassical shock waves (for further discussion regarding
the need for pointwise estimates, see \cite{HZ}).
It seems an observation (though elementary) of possibly wider 
use that our $H^2$ center stable manifold construction can be used to obtain
pointwise control on the solution in this way, extending a bit further
the basic ideas of \cite{Z5}.
As our study of conditional stability
was partly motivated by discussions in 
\cite{GZ,AMPZ,Z7} of certain unstable undercompressive shocks
and their effect on solution structure through metastable behavior,
it seems desirable to fit such nonclassical waves in the
theory.

An interesting open problem is to investigate conditional stability
of a planar standing shock $u(x,t)\equiv \bar u(x_1)$ of a
multidimensional system of conservation laws
$$
u_t + \sum_j f_j(u)_{x_j}= \sum_{jk} (b_{jk}(u)u_{x_k})_{x_j}
$$
which likewise (by the multidimensional arguments of \cite{Z1,Z2,Z3})
reduces to construction of a center stable manifold, in this case
involving an infinite-dimensional unstable subspace corresponding
to essential spectra of the linearized operator $L$ about the wave.

\medskip

{\bf Plan of the paper.}
In Section \ref{s:existence}, 
we establish existence of center stable manifolds
for general quasilinear parabolic PDE.
In Section \ref{s:cond}, we establish conditional stability
on the center stable manifold by a modification of the pointwise arguments 
of \cite{HZ,RZ} in the stable ($p=0$) case.


\section{Existence of Center Stable Manifold}\label{s:existence}

Defining the perturbation variable $v:=u-\bar u$, we
obtain after a brief computation the nonlinear perturbation
equations
\be\label{npert}
v_t - Lv=N(v),
\ee
where
\be\label{L}
Lv:= b(\bar u)v_{xx} + (db(\bar u)v) \bar u_x
-h_u(\bar u, \bar u_x)v
- h_{u_x}(\bar u, \bar u_x)v_x
\ee
denotes the linearized operator about the wave and
\ba\label{N1}
N(v)&:= 
\Big( b(\bar u+v)(\bar u+v)_{xx}-b(\bar u)\bar u_{xx} 
-b(\bar u)v_{xx}- (db(\bar u)v) \bar u_{xx} \Big) \\
&\quad
-\Big( h(\bar u+v,\bar u_x +v_x) - h(\bar u,\bar u_x) -h_u(\bar u, \bar u_x)v
- h_{u_x}(\bar u, \bar u_x)v_x\Big)\\
\ea
is a quadratic order residual.
We seek to construct a $C^k(H^2)$ local center stable manifold about the
equilibrium $v\equiv 0$, that is, a locally invariant $C^k$ manifold tangent
(in Frechet sense) with respect to $H^2$ to the center stable
subspace $\Sigma_{cs}$. 

\subsection{Preliminary estimates}

Denote by $\Sigma_u$ and $\Sigma_{cs}$ the unstable and center stable
subspaces of $L$ and $\Pi_u$ and $\Pi_{cs}$ the associated spectral projections.

\begin{proposition}[\cite{TZ1}]\label{linest}
Under assumptions (A0)--(A3), $L$ generates an analytic semigroup
$e^{Lt}$ satisfying
 \begin{equation} \label{bound-pde} 
\begin{aligned}
| e^{t L}\Pi_{cs}  |_{L^2\to L^2}  &\leq  C_\o   e^{\o t}, \\
| e^{-t L}\Pi_u  |_{L^2\to H^2}  &\leq  C_\o e^{-\b t},\\
\end{aligned}
 \end{equation}
 for some $\b > 0,$ and for all $\o > 0,$ for all $t\ge 0$.
\end{proposition}

\begin{proof}
Standard semigroup estimates for second-order 
elliptic operators; see \cite{TZ1,Z5}, or Appendix \ref{extraproofs}.
\end{proof}

Introducing a $C^\infty$ cutoff function
$
\rho(x) =\begin{cases}
1  & | x | \leq 1, \\ 
0 &  | x | \geq 2,\end{cases}
$
 let 
$$
 N^\delta (v) := 
\rho\Big( \frac{ | v |_{H^2}}{\delta}\Big) N(v).
$$

 \begin{lemma}[\cite{TZ1}] \label{trunc} 
Assuming (A0)--(A3), the map $N^\delta: H^2 
\to L^2 $ is $C^{k+1}$ and its Lipschitz norm with respect to $v$ 
is $O(\delta)$ as $\delta\to 0.$
Moreover, 
\be\label{Nquad}
|N^\delta(v)|_{H^2}\le C|v|_{H^2}^2.
\ee
 \end{lemma}

 \begin{proof}
See Appendix \ref{extraproofs}.
 \end{proof}

\begin{corollary}\label{integrand}
Under assumptions (A0)--(A3),
 \begin{equation} \label{bound-integrand} 
\begin{aligned}
| e^{t L}\Pi_{cs} N^\delta |_{H^2\to L^2}  &\leq  
C_\o e^{\o t}, \\
| e^{-t L}\Pi_u N^\delta |_{H^2\to H^2}  &\leq  C_\o e^{-\b t}, \\
\end{aligned}
 \end{equation}
for some $\b > 0,$ and for all $\o > 0,$ for all $t\ge 0$,
with Lipschitz bounds
 \begin{equation} \label{bound-integrand-lip} 
\begin{aligned}
| e^{t L}\Pi_{cs} dN^\delta |_{H^2\to L^2}  &\leq  
C_\o \delta  e^{\o t}, \\
| e^{-t L}\Pi_u dN^\delta |_{H^2\to H^2}  &\leq  C_\o  \delta e^{-\b t}. \\
\end{aligned}
 \end{equation}
\end{corollary}

\subsection{Fixed-point iteration scheme}\label{it}

Applying projections $\Pi_j$, $j=cs,u$ to the truncated equation
\be\label{qpert}
v_t-Lv=N^\delta(v),
\ee
we obtain using the variation of constants formula equations
$$
\Pi_j v(t)=
e^{L(t-t_{0,j})}\Pi_j v(t_{0,j}) + \int_{t_{0,j}}^t e^{L(t-s)}\Pi_j  N^\delta(v(s))\,ds,
$$
$j=cs,u$, so long as the solution $v$ exists,
with $t_{0,j}$ arbitrary.
Assuming growth of at most $|v(t)|_{H^2}\le Ce^{\tilde \theta t}$ in positive 
time, we find for $j=u$ using bounds \eqref{bound-pde}(ii)
and \eqref{bound-integrand-lip}(ii) 
that, as $t_{0,u}\to +\infty$, the first term
$e^{L(t-t_{0,u})}\Pi_u v(t_{0,u})$ converges to zero while
the second, integral term
converges to $\int_{t}^{+\infty} e^{L(t-s)}\Pi_u  N^\delta(v(s))\,ds$,
so that, denoting $w:=\Pi_{cs}v$, $z:=\Pi_u v$, we have
\be\label{vu}
z(t)=\cT(z,w)(t):=  -\int_t^{+\infty} e^{L(t-s)}\Pi_u  N^\delta( (w+z)(s))\,ds.
\ee

Likewise, choosing $t_{0,cs}=0$, we have
\be\label{vcs}
w(t)=
e^{Lt}\Pi_{cs} w_0 + \int_{0}^t e^{L(t-s)}\Pi_{cs}  N^\delta( (w+z)(s))\,ds,
\ee
$w_0:=\Pi_{cs}v(0)$.
On the other hand, 
we find from the original differential equation projected onto
the center stable component, after some rearrangement, 
that $w$ satisfies the Cauchy problem
\ba\label{diffvcs}
w_t- b(\bar u)w_{xx}&=\Pi_{cs}b(\bar u)z_{xx}- \Pi_u b(\bar u)w_{xx}\\
&\quad
+ \Pi_{cs}M(w+z)
+\Pi_{cs}N^\delta(w+z)\\
\ea
with initial data $ w_0=\Pi_{cs}v(0)$ given at $t=0$, where
\be\label{M}
M(v):= (db(\bar u)v) \bar u_x -h_u(\bar u, \bar u_x)v 
- h_{u_x}(\bar u, \bar u_x)v_x.
\ee
We shall use these two representations together to obtain optimal
estimates, the first for decay, through standard linear semigroup
estimates, and the second for regularity,
through the nonlinear damping estimate \eqref{damp1} below.

Viewing \eqref{vcs}, or alternatively \eqref{diffvcs},
as determining $w=\cW(z,w_0)$ as a function of $z$, we seek $z$ as
a solution of the fixed-point equation
\be\label{fixedpt}
z=\tilde \cT(z,w_0):=\cT(z, \cW(z,w_0)).
\ee
As compared to the standard ODE construction of, e.g., \cite{B,VI,TZ1,Z5},
in which \eqref{vu}--\eqref{vcs} together are considered as a fixed-point
equation for the joint variable $(w,z)$, this amounts to treating
$w$ {\it implicitly}.
This is a standard device in situations of limited regularity;
see, e.g., \cite{CP,GMWZ,RZ}.

It remains to show, first, that $\cW$, hence $\tilde\cT$, is well-defined 
on a space of slowly-exponentially-growing functions and,
second, that $\tilde\cT$ is contractive on that space, 
determining a $C^k$ solution
$z=z(w_0)$ similarly as in the usual ODE construction.
We carry out these steps in the following subsections.

\subsection{Nonlinear energy estimates}\label{nde}

\bl[\cite{Z5}]\label{projlem1}
Under assumptions (A0)--(A3), 
for all $1\le p\le \infty$, $0\le r\le 4$, 
\ba
|\Pi_{u}|_{L^p\to W^{r,p}} \, , \;
|\Pi_{cs}|_{W{r,p}\to W^{r,p}}&\le C.
\ea
\el

\begin{proof} See Appendix \ref{extraproofs}
\end{proof}

\bpr\label{damp1}
Under assumptions (A0)--(A3), for $\delta$ sufficiently small,
if the solution of \eqref{diffvcs} exists on $t\in [0,T]$,
then, for some constants $\theta_{1,2}>0$, 
and all $0\leq t\leq T$,
\begin{equation}\label{Ebounds1}
|w(t)|_{H^2}^2 \leq C e^{-\theta_1 t} |w_0|^2_{H^2} 
+ C \int_0^t e^{-\theta_2(t-s)} (|w|_{L^2}^2 + |z|_{L^2}^2) (s)\,ds.
\end{equation}
\epr

\begin{proof}
Let us first consider the simpler case that $b$ is uniformly elliptic,
$$
\Re b:= \frac{1}{2}(b+ b^t)\ge \theta >0.
$$
Taking the $L^2$ inner product 
in $x$ of $\sum_{j=0}^2\partial_x^{2j}w$ against \eqref{diffvcs}, 
integrating by parts and rearranging the resulting terms, 
we obtain 
$$
\begin{aligned}
\partial_t |w|_{H^2}^2(t) 
& \leq -\langle \partial_x^3 w, b(\bar u)\partial_x^3 w \rangle
+ C \big(|w|_{H^2}^2 + |z|_{H^4}^2\big),\\
& \leq -\theta |\partial_x^{3} w|_{L^2}^2 
+ C \big(|w|_{H^2}^2 + |z|_{H^4}^2\big),\\
\end{aligned}
$$
$\theta>0$, for $C>0$ sufficiently large, where we have repeatedly
used the bounds of Lemma \ref{trunc} and of 
Lemma \ref{projlem1} with $p=2$ to absorb the error terms coming
from the righthand side of \eqref{diffvcs}, and have used Moser's inequality
to bound $|\partial_x N^\delta(v)|_{L^2}\le C|v|_{L^\infty}|v|_{H^3}
\le  C\delta |v|_{H^3}$ whenever $N^\delta$ does not vanish,
so that $|v|_{H^2}\le 2\delta$.
Using the Sobolev interpolation
\[
|w|_{H^2}^2 \leq \tilde{C}^{-1} |\partial_x^{3} w|_{L^2}^2 + \tilde{C} | w|_{L^2}^2
\]
for $\tilde{C}>0$ sufficiently large, and
observing that
\be\label{eqnorm}
 |z|_{H^4}\le C|z|_{L^2} 
\ee
by equivalence of norms on finite-dimensional spaces, we obtain 
\[
\partial_t |w|_{H^2}^2(t) \leq -\tilde{\theta} |w|_{H^2}^2 + C \left(|w|_{L^2}^2 + |z|_{L^2}^2\right),
\]
from which \eqref{Ebounds1} follows by Gronwall's inequality.

To treat the general case, we note that
$\Re \sigma(b)\ge \theta>0$ by Lyapunov's Lemma implies that there exists
a positive definite matrix $P$ such that $\Re(Pb)\ge \frac{2\theta}{3}$,
whence, by a partition of unity argument, there exists a smooth positive
definite matrix-valued function $P(u)$ such that $\Re (Pb)\ge \frac{\theta}{2}$.
Taking the $L^2$ inner product 
in $x$ of $\sum_{j=0}^2\partial_x^{j}P(\bar u)\partial_x^jw$ 
against \eqref{diffvcs}, 
integrating by parts and rearranging the resulting terms, 
we obtain
\begin{equation}\label{en1}
\begin{aligned}
\partial_t \cE(w)
& \leq -\langle \partial_x^3 w, Pb(\bar u)\partial_x^3 w \rangle
+ C \big(|w|_{H^2}^2 + |z|_{H^4}^2\big),\\
& \leq -\frac{\theta}{2} |\partial_x^{3} w|_{L^2}^2 
+ C \big(|w|_{H^2}^2 + |z|_{H^4}^2\big),\\
\end{aligned}
\end{equation}
$\theta>0$, where
$\cE(w):=\sum_{j=0}^2 \langle \partial_x^j w, 
Pb(\bar u) \partial_x^j w \rangle$ is equivalent to $|\cdot|_{H^2}^2$.
By Sobolev interpolation and \eqref{eqnorm}, we have therefore
$$
\begin{aligned}
\partial_t \cE(w)
& \leq -{\tilde \theta}\cE
+ C \big(|w|_{L^2}^2 + |z|_{L^2}^2\big),\\
\end{aligned}
$$
$\tilde \theta>0$, from which \eqref{Ebounds1} follows 
again by Gronwall's inequality.
\end{proof}

\bpr\label{Lip1}
Under assumptions (A0)--(A3),
if solutions $w_1,z_1$ and $w_2,z_2$ of \eqref{diffvcs} exist for all
$t\ge 0$,
then, for $\delta$ sufficiently small, any constant $\theta \ge 0$ and some $C=C(\theta)$, 
for all $t\ge 0$,
\begin{equation}\label{LEbounds1}
\begin{aligned}
 \int_0^t e^{-\theta s} & |w_1-w_2|_{H^2}^2(s) \,ds \le
C|w_{0,1}-w_{0,2}|_{H^1}^2\\
&\quad +
C \int_0^t e^{-\theta s} (|w_1-w_2|_{L^2}^2 + |z_1-z_2|_{L^2}^2) (s)\,ds ,
\end{aligned}
\end{equation}
\begin{equation}\label{LEbounds2}
\begin{aligned}
 \int_t^{+\infty} e^{\theta (t-s) } & |w_1-w_2|_{H^2}^2(s) \,ds \le
C|w_1-w_2|_{H^1}^2(t)\\
&\quad +
C \int_t^{+\infty} e^{\theta (t- s)} 
(|w_1-w_2|_{L^2}^2 + |z_1-z_2|_{L^2}^2) (s)\,ds \\
&
\qquad \qquad \qquad
\quad
\le
Ce^{-\theta t}|w_{0,1}-w_{0,2}|_{H^1}^2\\
&\quad 
+C e^{2\eta t}(\|w_1-w_2\|_{L^2_{-\eta}}^2 + \|z_1-z_2\|_{L^2_{-\eta}}^2 ).
\end{aligned}
\end{equation}
\epr

\begin{proof}
Subtracting the equations for $w_1,z_1$ and $w_2,z_2$, we obtain,
denoting $\dot w:=w_1-w_2$, $\dot z:=z_1-z_2$, the equation
\ba\label{vardiffvcs}
\dot w_t- b(\bar u)\dot w_{xx}&=
\Pi_{cs}b(\bar u)\dot z_{xx}- \Pi_u b(\bar u)\dot w_{xx}\\
&\quad
+ \Pi_{cs}M(\dot w+\dot z)
+\Pi_{cs}\big( N^\delta(w_1+z_1) -N^\delta(w_2+z_2) \big)\\
\ea
with initial data $ \dot w_0=w_{0,1}-w_{0,2}$ at $t=0$.

Performing an $H^1$ version of the energy estimate in the 
proof of Proposition \ref{damp1}-- that is,
taking the $L^2$ inner product 
in $x$ of $\sum_{j=0}^1\partial_x^{j}P(\bar u)\partial_x^j\dot w$ 
against \eqref{vardiffvcs}, $P$ as in the proof of Proposition \ref{damp1},
integrating by parts, and rearranging the resulting terms--
we obtain
\begin{equation}\label{varen1}
\begin{aligned}
\partial_t \cE(\dot w)
& \leq -\frac{\alpha}{2} |\partial_x^{2} \dot w|_{L^2}^2 
+ C \big(|\dot w|_{H^1}^2 + |\dot z|_{H^2}^2\big),\\
\end{aligned}
\end{equation}
$\alpha>0$, where
$\cE(\dot w):=\sum_{j=0}^1 \langle \partial_x^j \dot w, 
Pb(\bar u) \partial_x^j \dot w \rangle$ is equivalent to $|\dot w|_{H^1}^2$.
Here, we have used in a key way the $H^2\to L^2$ Lipschitz bound on $N^\delta$
to bound 
$$
| N^\delta(w_1+z_1) -N^\delta(w_2+z_2) |_{L^2}\le 
C\delta(|\dot w|_{H^2}+|\dot z|_{H^2}),
\qquad \delta <<1,
$$
and thus
$
\langle
|\partial_x \dot w, Pb(\bar u)\partial_x \dot w \rangle ,
N^\delta(w_1+z_1) -N^\delta(w_2+z_2) \rangle|
\le C(|\dot w|_{H^2}^2+|\dot z|_{H^2}^2).
$

By Sobolev interpolation and \eqref{eqnorm}, we have therefore
\ba\label{sob1}
\partial_t \cE(\dot w) + |\dot w|_{H^2}^2
& \leq 
 C \big(|\dot w|_{L^2}^2 + |\dot z|_{L^2}^2\big),\\
\ea
whence \eqref{LEbounds1} follows 
by Gronwall's inequality, $\cE\ge 0$, and
$| \cE(\dot w_0)|\le C|\dot w_0|_{H^1}^2.  $

The first line of \eqref{LEbounds2} follows similarly, using
$| \cE(\dot w(t))|\le C|\dot w|_{H^1}^2(t).  $
Noting that we could in place of \eqref{sob1} have rearranged 
\eqref{varen1} as
$$
\begin{aligned}
\partial_t \cE(\dot w) 
& \leq  -\frac{\alpha}{2} \cE(\dot w)+
 C \big(|\dot w|_{L^2}^2 + |\dot z|_{L^2}^2\big),\\
\end{aligned}
$$
and argued as in Proposition \ref{damp1} to obtain
\ba\label{H1LEbounds1}
|\dot w(t)|_{H^1}^2 &\leq C e^{-\theta_1 t} |\dot w_0|^2_{H^1} 
+ C \int_0^t e^{-\theta_2(t-s)} (|\dot w|_{L^2}^2 + |\dot z|_{L^2}^2) (s)\,ds\\
&\le
C e^{-\theta_1 t} |\dot w_0|^2_{H^1} 
+C e^{2\eta t}(\|\dot w\|_{L^2_{-\eta}} + \|\dot z\|_{L^2_{-\eta}} ),
\ea
substituting in the first line
of \eqref{LEbounds2}, and estimating similarly
$$
 \int_t^{+\infty} e^{\theta (t- s)} 
(|\dot w|_{L^2}^2 + |\dot z|_{L^2}^2) (s)\,ds 
\le
C e^{2\eta t}(\|\dot w\|_{L^2_{-\eta}} + \|\dot z\|_{L^2_{-\eta}} ),
$$
we obtain the second line of \eqref{LEbounds2}, completing the proof.
\end{proof}

\br
\textup{
The absence of a uniform $H^3\to H^1$ Lipschitz bound on $N^\delta$
prevents us from obtaining a pointwise $H^2$ energy estimate on $\dot w$
like the one obtained on $w$ in Proposition \ref{damp1}.
}
\er

\subsection{Basic existence result}\label{pdecsproof}
Define now the negatively-weighted sup norm 
$$
\|f\|_{-\eta}:=\sup_{t\ge 0}e^{-\eta t}|f(t)|_{H^2},
$$
noting that $ |f(t)|_{H^2}\le e^{\tilde \theta t} \|f\|_{-\tilde \theta}$ 
for all $t\ge 0$, and denote by $\cB_{-\eta}$ the Banach space
of functions bounded in $\|\cdot\|_{-\eta}$ norm.
Define also the auxiliary norm 
$$
\|f\|_{L^2_{-\eta}}:=\sup_{t\ge 0}e^{-\eta t}|f(t)|_{L^2}.
$$

\bl\label{Texist}
Under assumptions (A0)--(A3),
for $3\o<\eta<\beta$ and $\delta>0$ and $w_0\in H^2$ sufficiently small,
for each $z\in \cB_{-\eta}$,
there exists a unique solution $w=:\cW(z,w_0)\in \cB_{-\eta}$ 
of \eqref{vcs}, \eqref{diffvcs}, with
\be\label{comp}
\|w\|_{-\eta} \le C(|w_0|_{H^2}+\delta\|z\|_{-\eta})
\ee
and 
\be\label{Lipcomp}
\|\cW(z_1, w_{0,1})-\cW(z_2, w_{0,2})\|_{L^2_{-\eta}}
\le C\delta \|z_1-z_2\|_{-\eta}
+ C\|w_{0,1}-w_{0,2}\|_{-\eta}.
\ee
\el

\begin{proof}
{\it (i) ($H^2$ bound})
By short-time $H^2$-existence theory, and the earlier-observed fact
\eqref{eqnorm}, an $H^2$ solution of \eqref{diffvcs} exists and
remains bounded $H^2$ up to some time $T>0$ provided
that $w_0$ is bounded in $H^2$, whereupon \eqref{Ebounds1}
holds.

Using now the integral representation \eqref{vcs}, and applying
\eqref{bound-pde}(i), \eqref{bound-integrand-lip}(i)
$\Rightarrow |N^\delta (v(t))|\le \delta |v(t)|$, and
$
|w(t)|_{H^2}\le e^{\eta |t|} \|w\|_{-\eta,T}:=
\sup_{0\le s\le T}e^{-\eta s}|w(s)|_{H^2},
$
we obtain for $0\le t\le T$ that
\ba\nonumber
|w(t)|_{L^2}&\le C e^{\o |t|}|w_{cs}|_{L^2}
+C\delta (\|w\|_{-\eta,T }+ \|z\|_{-\eta })
\int_{0}^t e^{  \o |t-s|} e^{\eta |s|}\, ds,
\ea
hence, using $\o\pm \eta >0$ that
$
|w(t)|_{L^2}\le C e^{\eta |t|}
\Big(
|w_{cs}|_{L^2} +\delta (\|w\|_{-\eta,T }+ \|z\|_{-\eta }) \Big).
$

Applying \eqref{Ebounds1}, we then obtain
\be\nonumber
|w(t)|_{H^2}\le C e^{\eta |t|}
\Big(
|w_{cs}|_{H^2} +\delta (\|w\|_{-\eta,T }+ \|z\|_{-\eta }) \Big)
\ee
for all $0\le t\le T$, and thus
$
\|w\|_{-\eta,T}\le 
C\Big( |w_{cs}|_{H^2} +\delta (\|w\|_{-\eta,T }+ \|z\|_{-\eta }) \Big),
$
whence, for $\delta$ sufficiently small,
$\|w\|_{-\eta,T}\le C( |w_{cs}|_{H^2} + \|z\|_{-\eta })$. 
Since this bound is independent of $T$, we obtain by continuation
global existence of $w$ and, letting $T\to \infty$, \eqref{comp}
as claimed.

{\it (ii) ($H^2\to L^2$ Lipschitz bounds)}
Now consider a pair of data $z_1, w_{0,1}$ and
$z_2, w_{0,2}$, and compare the resulting solutions,
denoting 
$
(\dot z, \dot w , \dot w_0):=(z_1-z_2,w_1-w_2,w_{0,1}-w_{0,2}).
$
Using the integral representation \eqref{vcs}, and applying
\eqref{bound-pde}(i), \eqref{bound-integrand-lip}(i),
and the definition of $\|\cdot\|_{-\eta}$, we obtain
for all $t\ge 0$ 
\ba\label{tempz}
|\dot w(t)|_{L^2}&\le C e^{\o |t|}|\dot w_{0}|_{L^2}
+ C\delta \int_{0}^t e^{  \o |t-s|} 
 (|\dot w\dot |_{H^2 }+ |\dot z|_{H^2 })(s) \, ds.
\ea

Estimating  as before
$$
 C\delta \int_{0}^t e^{  \o (t-s)} |\dot z|_{H^2 }(s) \, ds
\le
 C\delta  \|\dot z\|_{-\eta}\int_{0}^t
e^{  \o (t-s)} e^{\eta s} \, ds
\le
 C\delta  \|\dot z\|_{-\eta} e^{\eta t}
$$
and, by \eqref{LEbounds1} with $\theta=3\o$ together with the
Cauchy--Schwarz inequality, 
$$ 
\begin{aligned}
C\delta \int_{0}^t e^{  \o (t-s)} 
 |\dot w\dot |_{H^2 }(s) \, ds &\le 
C\delta \Big(\int_{0}^t e^{  -\o (t-s)} \Big)^{1/2}
\Big(\int_{0}^t e^{3  \o (t-s)} 
 |\dot w\dot |_{H^2 }^2(s) \, ds \Big)^{1/2}\\
 &\le 
C\delta \Big(e^{3\o t}|w_0|_{H^1}^2+
\int_{0}^t e^{ 3 \o (t-s)} 
 (|\dot w\dot |_{L^2 }^2+ |\dot z|_{L^2}^2)(s) \, ds \Big)^{1/2},
\end{aligned}
$$
we obtain, substituting in \eqref{tempz}, 
\ba\label{tempz2}
|\dot w(t)|_{L^2}^2&\le C e^{6\o |t|}|\dot w_{0}|_{H^1}^2
+ C\delta^2 \Big(  \|\dot z\|_{-\eta}^2 e^{2\eta t}
+  \int_{0}^t e^{  3\o (t-s)} 
 (|\dot w\dot |_{L^2 }^2+ |\dot z\dot |_{L^2 }^2) (s) \, ds \Big)\\
&\le
 e^{2\eta t} \Big( C|\dot w_{0}|_{H^1}^2
+ C\delta^2 ( \|\dot w\|_{L^2_{-\eta}}^2 + \|\dot z\|_{-\eta}^2)\Big).
\ea

This yields 
\ba\nonumber
\|\dot w(t)\|_{L^2_{-\eta}}&\le C |\dot w_{0}|_{-\eta}
+ C\delta ( \|\dot w\|_{L^2_{-\eta}} + \|\dot z\|_{-\eta}),
\ea
from which \eqref{Lipcomp} follows by smallness of $\delta$.
\end{proof}

\begin{proof}[Proof of Proposition \ref{t:cspde}.] 
{\it (i) (Boundedness on a ball)}
Again, working in $\cB_{-\eta}$, recall that
\be\nonumber
z(t)=\tilde \cT(z,w_0)(t):=  -\int_t^{+\infty} e^{L(t-s)}\Pi_u  N^\delta( (w+z)(s))\,ds,
\ee
where $w=\cW(z,w_0)$.
Using \eqref{vu}, and applying
\eqref{bound-pde}(ii), \eqref{bound-integrand-lip}(ii),
$|N^\delta (v(t))|\le \delta |v(t)|$, and
$ |z(t)|_{H^2} \le e^{\eta |t|} \|z\|_{-\eta} $
and, by \eqref{comp},
$$
|w(t)|_{H^2}
\le e^{\eta |t|} \|w\|_{-\eta}\le C e^{\eta |t|} 
(|w_0|_{H^2}+\delta\|z\|_{-\eta}),
$$
we obtain for $0\le t$ that
\ba\nonumber
|\tilde T(z,w_0)(t)|&\le C\delta 
(|w_0|_{H^2}+\delta\|z\|_{-\eta})
\int_t^{+\infty} e^{\beta(t-s)} e^{\eta |s|}\, ds,\\
\ea
hence, using $\beta\pm \eta >0$
and taking $\delta $ and $|w_0|_{H^2}$ sufficiently small, that
$\tilde \cT(\cdot, w_0)$ maps the ball $B(0,r)$ to itself, for
$r>0$ arbitrarily small but fixed.

{\it (ii) (Lipschitz bounds)}
Similarly, we find using \eqref{Lipcomp}, \eqref{bound-integrand-lip}(ii), 
and Lemma \ref{trunc} that
\ba\label{tempz5}
\|\tilde \cT(z_1,w_{0,1})&-\tilde \cT(z_2,w_{0,2})\|_{-\eta}\\
&
\le \sup_t C\delta  e^{-\eta t}\int_t^{+\infty} e^{\beta(t-s)} (|\dot w|_{H^2}+|\dot z|_{H^2})(s) \, ds
\\ &
\le \sup_t C_1 \delta
\Big( \|\dot z\|_{-\eta }
+
e^{-\eta t} \int_t^{+\infty} e^{\beta(t-s)} |\dot w|_{H^2}(s)  ds \Big).
\\
\ea
Using the Cauchy--Schwarz inequality and
\eqref{LEbounds2} with $\theta=\beta$ to estimate
\ba\nonumber
\int_t^{+\infty} e^{\beta(t-s)} |\dot w|_{H^2}(s)  ds &\le
\Big(  \int_t^{+\infty} e^{\beta(t-s)} ds \Big)^{1/2}
\Big(  \int_t^{+\infty} e^{\beta(t-s)} |\dot w|_{H^2}^2(s)  ds \Big)^{1/2}
\\
&\le
C_3\Big(
|\dot w|_{H^1}^2(t) +
 \int_0^t e^{\beta (t-s)} (|\dot w|_{L^2}^2 + |\dot z|_{L^2}^2) (s)ds
\Big)^{1/2}\\
&\le
C_3e^{\eta t}\Big(|\dot w_{0}|_{H^1}
+\|\dot w\|_{L^2_{-\eta}} + \|\dot z\|_{L^2_{-\eta}} \Big)
\\
\ea
and applying \eqref{Lipcomp}, we obtain 
\ba\label{finalLip}
\|\tilde \cT(z_1,w_{0,1})-\tilde \cT(z_2,w_{0,2})\|_{-\eta}&\le
C\delta ( \|\dot w_0\|_{-\eta} + \|\dot z\|_{L^2_{-\eta}}).
\ea
This yields at once contractivity on $B(0,r)$,
hence existence of a unique fixed point $z=\cZ(w_0)$, and
Lipshitz continuity of $\cZ$ from $\Sigma_{cs}$ to $\cB_{-\eta}$,
by the Banach Fixed-Point Theorem
with Lipschitz dependence on parameter $w_0$.

{\it (ii) (Existence of a Lipschitz invariant manifold)}
Defining
\ba\label{Phidefc}
\Phi(w_{0})&:= \cZ(w_0)|_{t=0} =
 -\int_{0}^{+\infty}e^{-Ls}\Pi_{u}  N^\delta (v(s))\,ds,
\ea
we obtain a Lipschitz function from $\Sigma_{cs}\to \Sigma_{u}$,
whose graph over $B(0,r)$ is the invariant manifold of solutions 
of \eqref{qpert} growing at exponential rate $|v(t)|\le Ce^{\eta t}$ 
in forward time.  From the latter characterization, we obtain evidently
invariance in forward time.
Since the truncated equations \eqref{qpert} agree with the original
PDE so long as solutions remain small in $H^2$, this gives local invariance
with respect to \eqref{soln} as well.
By uniqueness of fixed point solutions, we have $\cZ(0)=0$ and
thus $\Phi_{cs}(0)=0$, so that the invariant manifold passes through
the origin.
Likewise, any bounded, sufficiently small solution of \eqref{soln} 
in $H^2$ is a bounded small soluton of \eqref{qpert} as well, so
by uniqueness is contained in the center stable manifold.

{\it (ii) (Quadratic-order tangency)}
By \eqref{Phidefc}, \eqref{bound-pde}, \eqref{Nquad}, 
and \eqref{comp},
\ba\nonumber
|\Phi(w_{cs})|_{H^2} &=
 \Big|\int_{0}^{+\infty}e^{-Ls}\Pi_{u}  N^\delta (w+z)(s)\,ds
\Big|_{H^2}\\
&\le C\int_0^{+\infty} e^{-\beta s}(|\cW|_{H^2}^2 +|\cZ|_{H^2}^2)(s)ds\\
&\le C_1(\|\cW\|_{-\eta}^2 +\|\cZ\|_{-\eta}^2)
\le C_2(|w_{cs}|_{H^2}^2 +\|\cZ\|_{-\eta}^2).\\
\ea
By $\cZ(0)=0$ and Lipshitz continuity of $\cZ$, we have
$\|\cZ\|_{-\eta}\le C|w_{cs}|_{H^2}$, whence
\be\label{tanbd}
|\Phi(w_{cs})|_{H^2}\le C_3 |w_{cs}|_{H^2}^2,
\ee
verifying quadratic-order tangency at the origin.
\end{proof}

\subsection{Translation-invariance}\label{s:trans}

We conclude by indicating briefly how to recover translation-invariance
of the center stable manifold, following \cite{Z5,TZ1}.
Differentiating with respect to $\alpha$ the relation
$ 0\equiv = \cF(\bar u(x+\alpha)) $,
we recover the standard fact that $ \phi := \bar u_x $
is an $L^2$ zero eigenfunction of $L$, by the assumed decay of $\bar u_x$.

Define orthogonal projections
\begin{equation}
\label{proj}
\Pi_2 :=
\frac{ \phi \, \langle \phi, \cdot\rangle}{|\phi|_{L^2}^2}, 
\qquad \Pi_1 := \Id - \Pi_2,
\end{equation}
onto the range of right zero-eigenfunction $\phi:=(\partial/\partial x)\bar u$
of $L$ and its orthogonal complement $\phi^\perp$ in $L^2$, where 
$\langle \cdot, \cdot \rangle$ denotes standard $L^2$ inner product.

\begin{lemma}\label{Pi}
Under the assumed regularity $h\in C^{k+1}$, $k\ge 2$,
$\Pi_j$, $j=1,2$ are bounded as operators from $H^s$ to itself
for $0\le s \le k+2$.
\end{lemma}

\begin{proof}
Immediate, by the assumed decay of $\phi=\bar u_x$ and derivatives.
\end{proof}

\begin{proof}[Proof of Theorem \ref{t:cstrans}]
Introducing the shifted perturbation variable
\be\label{shift}
v(x,t):= u(x+\alpha(t),t)-\bar u(x)
\ee
we obtain the modified nonlinear perturbation equation
\begin{equation}
\label{epspert}
\d_t v= Lv + N(v) -\d_t \alpha (\phi+ \partial_x v),
\end{equation}
where $ L := \frac{\d {\cal F}}{\d u}( \bar u)$ and
$N$ as in \eqref{N} is a quadratic-order Taylor remainder.

Choosing $\d_t\alpha$ so as to cancel $\Pi_2$ of the righthand
side of \eqref{epspert}, we obtain the {\it reduced equations}
\be\label{redv}
\d_t v= \Pi_1( Lv + N(v))
\ee
and
\be\label{alphaeq}
\d_t \alpha =
\frac{\pi_2(Lv + N(v))}
{1+ \pi_2 (\d_x v)}
\ee
$v\in \phi^\perp$, 
$\pi_2 v:= \langle \tilde \phi, v\rangle |\phi|_{L^2}^2$,
of the same regularity as the original equations.

Here, we have implicitly chosen $\alpha(0)$ so that
$v(0)\in \phi^\perp$, or
$$
\langle \phi, u_0(x+\alpha)-\bar u(x)\rangle=0.
$$
Assuming that $u_0$ lies in a sufficiently small tube about
the set of translates of $\bar u$, or $u_0(x)=(\bar u+w)(x-\beta)$
with $|w|_{H^2}$ sufficiently small, this can be done in a unique
way such that $\tilde \alpha:=\alpha-\beta$ is small, as determined
implicitly by 
$$
0=\cG(w,\tilde \alpha):=
\langle \phi, \bar u(\cdot+\tilde \alpha)-\bar u(\cdot)\rangle
+\langle \phi, w(\cdot+\alpha) \rangle,
$$
an application of the Implicit Function Theorem noting that
$$
\d_{\tilde \alpha}\cG(0,0)=
 \langle \phi, \partial_x \bar u(\cdot)\rangle=|\phi|_{L^2}^2\ne0.
$$
With this choice, translation invariance under our construction is clear,
with translation corresponding to a constant shift in $\alpha$
that is preserved for all time.

Clearly, \eqref{alphaeq} is well-defined so long as 
$|\d_x v|_{L^\infty}\le C|v|_{H^2}$
remains small, hence we may solve the $v$ equation independently
of $\alpha$, determining $\alpha$-behavior afterward to determine
the full solution 
$$
u(x,t)= \bar u(x-\alpha(t))+v(x-\alpha(t),t).
$$
Moreover, it is easily seen (by the block-triangular structure
of $L$ with respect to this decomposition) that
the linear part $\Pi_1L=\Pi_1 L\Pi_1$ of the $v$-equation possesses
all spectrum of $L$ apart from the zero eigenvalue associated with 
eigenfunction $\phi$.
Thus, we have effectively projected out this zero-eigenfunction, and
with it the group symmetry of translation.

We may therefore construct the center stable manifold for the reduced
equation \eqref{redv}, automatically obtaining translation-invariance
when we extend to the full evolution using \eqref{alphaeq}.
See \cite{TZ1,Z5} for further details.
\end{proof}

\subsection{Application to viscous shock waves}

\begin{proof} [Proof of Proposition \ref{hA}]
Clearly, (H0) implies (A0) by the form \eqref{shockcase} of $h$,
while (H1) and (A1) are identical.
Plugging $u=\bar u(x)$ into \eqref{cons}, we obtain the 
standing-wave ODE $f(\bar u)_x=(b(\bar u)\bar u_x)_x$, or, integrating
from $-\infty$ to $x$, the first-order system
$$
b(\bar u)\bar u_x= f(\bar u)-f(u_-).
$$
Linearizing about the assumed critical points $u_\pm$ yields
linearized systems $w_t=df(u_\pm)w$, from which we see that
$u_\pm$ are nondegenerate rest points by (H2),
as a consequence of which (A3) follows by standard ODE theory \cite{Co}.

Finally, linearizing PDE \eqref{rcd} about the constant solutions
$u\equiv u_\pm$, we obtain $w_t=L_\pm w:= -df(u_\pm) w_x-b(u_\pm)w_{xx}$.
By Fourier transform, the limiting operators $L_\pm$ have
spectra $\lambda^\pm_j(k)=\sigma (-ikA^\pm-b_\pm k^2)$, where the Fourier
wave-number $k$ runs over all of $\R$; in particular, $L_\pm$
have spectra of nonpositive real part by (H3).
By a standard result of Henry \cite{He}, the essential spectrum
of $L$ lies to the left of the rightmost boundary of the spectra
of $L_\pm$, hence we may conclude that
the essential spectrum of $L$ is entirely nonpositive.
As the spectra of $L$ to the right of the essential spectrum
by sectoriality of $L$, consists of finitely
many discrete eigenvalues, this means that the spectra of $L$
with positive real part consists of $p$ unstable eigenvalues,
for some $p$, verifying (A2).
\end{proof}

\section{Conditional stability analysis}\label{s:cond}

From now on we specialize to the case
\eqref{shockcase}, \eqref{cons} of a system of viscous conservation laws
$ u_t= \cF(u)= (b(u)u_{x})_x -f(u)_x $, with associated linearized operator
$L=\partial{\cF}{u}(\bar u)$.

\subsection{Projector bounds}\label{projbds}
Let $\Pi_u$ denote the eigenprojection of $L$ onto its
unstable subspace $\Sigma_u$, and $\Pi_{cs}={\rm Id}- \Pi_u$
the eigenprojection onto its center stable subspace $\Sigma_{cs}$. 

\bl[\cite{Z5}]\label{projlem}
Assuming (H0)--(H1),
\be\label{comm}
\Pi_j \partial_x= \partial_x \tilde \Pi_j
\ee
for $j=u,\, cs$ and, for all $1\le p\le \infty$, $0\le r\le 4$, 
\ba\label{pbd}
|\Pi_{u}|_{L^{p}\to W^{r,p}}, |\tilde \Pi_{u}|_{L^{p}\to W^{r,p}}&\le C,\\
|\tilde \Pi_{cs}|_{W{r,p}\to W^{r,p} }, 
\; |\tilde \Pi_{cs}|_{W{r,p}\to W^{r,p}}&\le C,
\ea
and
\ba\label{wtpbd}
|(1+|x|^2)^{\frac{3}{4}}\Pi_{u}f|_{H^4}&\le |f|_{L^p},\\
|(1+|x|^2)^{\frac{3}{4}}\Pi_{u}f|_{H^4}&\le
|(1+|x|^2)^{\frac{3}{4}}f|_{H^4}.
\ea
Moreover,
\be\label{ptproj}
|\Pi_{cs}f(x)|\le Ce^{-\theta |x|}\sup_y |e^{-\theta |y|}f(y)|.
\ee
\el

\begin{proof}
Recalling that $L$ has at most finitely many unstable eigenvalues,
we find that $\Pi_u$ may be expressed as
$$
\Pi_u f= \sum_{j=1}^p \phi_j(x) \langle \tilde \phi_j, f\rangle,
$$
where $\phi_j$, $j=1, \dots p$ are generalized right eigenfunctions of
$L$ associated with unstable eigenvalues $\lambda_j$, 
satisfying the generalized eigenvalue equation $(L-\lambda_j)^{r_j}\phi_j=0$,
$r_j\ge 1$, and $\tilde \phi_j$
are generalized left eigenfunctions.
Noting that $L$ is divergence form, and that $\lambda_j\ne 0$,
we may integrate $(L-\lambda_j)^{r_j}\phi_j=0$ over $\R$ to
obtain $\lambda_j^{r_j}\int \phi_j dx=0$ and thus $\int\phi_j dx=0$.
Noting that $\phi_j$, $\tilde \phi_j$ and derivatives decay exponentially
by standard theory \cite{He,ZH,MaZ1}, we find that
$$
\phi_j= \partial_x \Phi_j
$$
with $\Phi_j$ and derivatives exponentially decaying, hence
$$
\tilde \Pi_u f=\sum_j \Phi_j \langle \partial_x \tilde \phi, f\rangle.
$$
Estimating 
$$
|\partial_x^j\Pi_u f|_{L^p}=|\sum_j \partial_x^j\phi_j \langle \tilde \phi_j f
\rangle|_{L^p}\le
\sum_j |\partial_x^j\phi_j|_{L^p} |\tilde \phi_j|_{L^q} |f|_{L^p}
\le C|f|_{L^p}
$$
for $1/p+1/q=1$
and similarly for $\partial_x^r \tilde \Pi_u f$, we obtain the claimed
bounds on $\Pi_u$ and $\tilde \Pi_u$, from which the bounds on
$\Pi_{cs}={\rm Id}-\Pi_u$ and
$\tilde \Pi_{cs}={\rm Id}-\tilde \Pi_u$ follow immediately.
Bounds \eqref{wtpbd} and \eqref{ptproj} follow similarly.
\end{proof}

\subsection{Linear estimates}

Let $G_{cs}(x,t;y):=\Pi_{cs}e^{Lt}\delta_y(x)$ denote
the Green kernel of the linearized solution operator on
the center stable subspace $\Sigma_{cs}$.
Then, we have the following detailed pointwise bounds
established in \cite{TZ2,MaZ1}.

\begin{proposition}[\cite{TZ2,MaZ1}]\label{greenbounds}
Assuming (H0)--(H5), (D1)--(D2), 
the center stable Green function may be decomposed as 
$G_{cs}=E+\tilde G$, where 
\begin{equation}\label{E}
E(x,t;y)=\sum_{j=1}^\ell 
\frac{\partial \bar u^\alpha(x)}{\partial \alpha_j}_{|\alpha=\des}e_j(y,t),
\end{equation}
\begin{equation}\label{e}
  e_j(y,t)=\sum_{a_k^{-}>0}
  \left(\textrm{errfn }\left(\frac{y+a_k^{-}(t+1)}{\sqrt{4\beta_k^{-}(t+1)}}\right)
  -\textrm{errfn }\left(\frac{y-a_k^{-}t}{\sqrt{4\beta_k^{-}(t+1)}}\right)\right)
  l_{jk}^{-}(y)
\end{equation}
for $y\le 0$ and symmetrically for $y\ge 0$, with
\begin{equation}\label{ljkbounds}
|l_{jk}^\pm|\le C, \qquad
|(\partial/\partial y)l_{jk}^\pm|\le C\gamma e^{-\eta |y|}, 
\end{equation}
and
\begin{equation}\label{Gbounds}
\begin{aligned}
|\d_x^s\tilde G(x,t;y)|&\le  Ce^{-\eta(|x-y|+t)} \\
&+
C(t^{-\frac{s}{2}}+e^{-\theta|x|})
\Big(
\sum_{k=1}^n 
t^{-\frac{1}{2}}
e^{-(x-y-a_k^{-} t)^2/Mt} e^{-\eta x^+} \\
&+
\sum_{a_k^{-} > 0, \, a_j^{-} < 0} 
\chi_{\{ |a_k^{-} t|\ge |y| \}}
t^{-1/2} e^{-(x-a_j^{-}(t-|y/a_k^{-}|))^2/Mt}
e^{-\eta x^+} \\
&+
\sum_{a_k^{-} > 0, \, a_j^{+}> 0} 
\chi_{\{ |a_k^{-} t|\ge |y| \}}
t^{-1/2} e^{-(x-a_j^{+} (t-|y/a_k^{-}|))^2/Mt}
e^{-\eta x^-}\Big), \\
\end{aligned}
\end{equation}
\begin{equation}\label{Gybounds}
\begin{aligned}
|\d_x^s\partial_y \tilde G(x,t;y)|&\le  Ce^{-\eta(|x-y|+t)}\\
& + 
Ct^{-\frac{1}{2}}
(t^{-\frac{s}{2}}+e^{-\theta|x|}+ \gamma e^{-\theta|y|})
\Big( \sum_{k=1}^n 
t^{-1/2}e^{-(x-y-a_k^{-} t)^2/Mt} e^{-\eta x^+} \\
&+
\sum_{a_k^{-} > 0, \, a_j^{-} < 0} 
\chi_{\{ |a_k^{-} t|\ge |y| \}}
t^{-1/2} e^{-(x-a_j^{-}(t-|y/a_k^{-}|))^2/Mt}
e^{-\eta x^+} \\
&+
\sum_{a_k^{-} > 0, \, a_j^{+}> 0} 
\chi_{\{ |a_k^{-} t|\ge |y| \}}
t^{-1/2} e^{-(x-a_j^{+} (t-|y/a_k^{-}|))^2/Mt}
e^{-\eta x^-}\Big) \\
\end{aligned}
\end{equation}
for $0\le s\le 1$,
for $y\le 0$ and symmetrically for $y\ge 0$,
for some $\eta$, $C$, $M>0$, where 
$a_j^\pm$ are as in \eqref{aj},  $\beta_k^\pm>0$,
$x^\pm$ denotes the positive/negative
part of $x$,  indicator function $\chi_{\{ |a_k^{-}t|\ge |y| \}}$ is 
$1$ for $|a_k^{-}t|\ge |y|$ and $0$ otherwise,
and $\gamma=1$ in the mixed or undercompressive
case and $0$ in the pure Lax or overcompressive case.
\end{proposition}

\begin{proof}
This follows from the observation \cite{TZ2} that the contour integral 
(inverse Laplace Transform) representation of $G_{cs}$ is exactly
that for the full Green kernel in the stable case $p=0$, and that
the resolvent kernel satisfies the same bounds.  Thus, we obtain
the stated bounds by the same argument used in \cite{MaZ1} to
bound the full Green kernel in the stable case.
See Appendix \ref{extraproofs} for further discussion.
\end{proof}

\bc[\cite{HZ}]\label{ebds}
Assuming (H0)--(H5), (D1)--(D2), 
\begin{equation}\label{ebounds}
\begin{aligned}
|e_j(y,t)|&\le C\sum_{a_k^->0}
  \left(\textrm{errfn }\left(\frac{y+a_k^{-}t}{\sqrt{4\beta_k^{-}(t+1)}}\right)
  -\textrm{errfn }\left(\frac{y-a_k^{-}t}{\sqrt{4\beta_k^{-}(t+1)}}\right)\right),\\
|e_j (y,t) &- e_j (y,+\infty)| \le C \textrm{errfn} (\frac{|y|-at}{M\sqrt{t}}), 
\quad \text{some }\, a>0 \\
|\partial_t  e_j(y,t)|&\le C t^{-1/2} \sum_{a_k^->0} e^{-|y+a_k^-t|^2/Mt},\\
|\partial_y  e_j(y,t)|&\le C t^{-1/2} \sum_{a_k^->0} e^{-|y+a_k^-t|^2/Mt}
+ C\gamma e^{-\eta|y|}\\
&\times
  \left(\textrm{errfn }\left(\frac{y+a_k^{-}t}{\sqrt{4\beta_k^{-}(t+1)}}\right)
  -\textrm{errfn }\left(\frac{y-a_k^{-}t}{\sqrt{4\beta_k^{-}(t+1)}}\right)\right),\\
|\partial_y e_j (y,t) &- \partial_y e_j(y,+\infty)|
 \le  C t^{-1/2} \sum_{a_k^->0} e^{-|y+a_k^-t|^2/Mt} \\
|\partial_{yt}  e_j(y,t)|&\le C
(t^{-1}+\gamma t^{-1/2}e^{-\eta|y|}) \sum_{a_k^->0} e^{-|y+a_k^-t|^2/Mt}\\
\end{aligned}
\end{equation}
for $y\le 0$, and symmetrically for $y\ge 0$, where $\gamma$ as above
is one for undercompressive profiles and zero otherwise. 
\ec

\begin{proof}
Direct computation using definition \eqref{e}; 
see \cite{MaZ1}. 
\end{proof}

\subsection{Convolution bounds}
From the above pointwise bounds, there follow by direct computation
the following convolution estimates established in \cite{HZ}, here
stated without proof.

\begin{lemma}[Linear 
estimates \cite{HZ}]\label{iniconvolutions}
Assuming (H0)--(H5), (D1)--(D2), 
\begin{equation}\label{iniconeq}
\begin{aligned}
\int_{-\infty}^{+\infty}|\tilde G(x,t;y)|(1+|y|)^{-3/2}\, dy
&\le C (\theta+\psi_1+\psi_2)(x,t),\\
\int_{-\infty}^{+\infty}|\tilde G_x(x,t;y)|(1+|y|)^{-3/2}\, dy
&\le C (t^{-\frac{1}{2}}+1) (\theta+\psi_1+\psi_2)(x,t),\\
\int_{-\infty}^{+\infty}|e_t(y,t)|(1+|y|)^{-3/2}\, dy
&\le C(1+t)^{-3/2},\\
\int_{-\infty}^{+\infty}|e(y,t)|(1+|y|)^{-3/2}\, dy
&\le C,\\
 \int^{+\infty}_{-\infty} |e(y,t)-e(y,+\infty)| (1+|y|)^{-3/2}\, dy 
&\le C(1+t)^{-1/2},\\
\end{aligned}
\end{equation}
for $0\le t\le +\infty$, $C>0$, where $\tilde G$ and 
$e$ are defined as in Proposition \ref{greenbounds}.
\end{lemma}

\begin{lemma}[Nonlinear 
estimates \cite{HZ}]\label{convolutions}
Under (H0)--(H5), (D1)--(D2), 
\begin{equation}\label{coneq}
\begin{aligned}
\int_0^t\int_{-\infty}^{+\infty}|\tilde G_y(x,t-s;y)|\Psi(y,s)\, dy ds
&\le C(\theta+\psi_1+\psi_2)(x,t),\\
\int_0^{t-1}\int_{-\infty}^{+\infty}|\tilde G_{yx}(x,t-s;y)|\Psi(y,s)\, dy ds
&\le C(\theta+\psi_1+\psi_2)(x,t),\\
\int_{t-1}^t\int_{-\infty}^{+\infty}|\tilde G_x(x,t-s;y)| 
(\theta+\psi_1+\psi_2)(y,s) \, dy ds
&\le C(\theta+\psi_1+\psi_2)(x,t),\\
\int_0^t\int_{-\infty}^{+\infty}|e_{yt}(y,t-s)|\Psi(y,s)\, dy ds
&\le C(1+t)^{-1},\\
\int_t^{+\infty} \int_{-\infty}^{+\infty}|e_y(y,+\infty)|\Psi(y,s)\, dy
&\le C\gamma(1+t)^{-1/2},\\
\int_0^t\int_{-\infty}^{+\infty}
|e_y(y,t-s)- e_y(y,+\infty)| \Psi(y,s)\, dyds
&\le C(1+t)^{-1/2},\\
\end{aligned}
\end{equation}
for $0\le t\le +\infty$, $C>0$, where $\tilde G$ and 
$e$ are as in Proposition \ref{greenbounds} and
\begin{equation}\label{source}
\begin{aligned}
\Psi(y,s)&:=
(1+s)^{1/2}s^{-1/2}(\theta + \psi_1+\psi_2)^2(y,s)\\
&\qquad +
(1+s)^{-1} (\theta+\psi_1+\psi_2)(y,s).\\
\end{aligned}
\end{equation}
\end{lemma}

We have by standard short-time theory the following additional bound.

\begin{lemma}[Commutator bound]\label{inicomm}
Assuming (H0)--(H5), (D1)--(D2), 
\begin{equation}\label{inicommeq}
\begin{aligned}
\int_{-\infty}^{+\infty}|(\tilde G_x+\tilde G_y)(x,t;y)|(1+|y|)^{-3/2}\, dy
&\le C (\theta+\psi_1+\psi_2)(x,t),\\
\ea
for $0\le t\le 1$, $C>0$, where $\tilde G$
is as in Proposition \ref{greenbounds}.
\end{lemma}

\begin{proof}
By standard short-time parametrix bounds \cite{Fr} on the
entire Green function, we have for $0\le t\le 1$ 
that $|G_x+G_y|\le Ct^{-\frac{1}{2}}
e^{-\frac{|x-y|^2}{Mt}}$ is of the same order order as $\tilde G$,
whence
\begin{equation}\nonumber
\begin{aligned}
\int_{-\infty}^{+\infty}|(G_x+G_y)(x,t;y)|(1+|y|)^{-3/2}\, dy
&\le C(1+|x|)^{-3/2}.\\
\ea
By direct computation, integration against
$|\partial_xG_u|$, $|\partial_y G_u|$, $|\partial_x E|$, or
$|\partial_y E|$ gives a contribution
that is also bounded by $ C(1+|x|)^{-3/2}$, yielding the result.
\end{proof}

\subsection{Reduced equations II}
We now restrict to the pure Lax or undercompressive case
$\ell=1$, following the simple stability argument of Section 3, \cite{HZ}.
The pure Lax or overcompressive case may be carried out following the
similar but slightly more complicated argument of Section 4, \cite{HZ}.

Define similarly as in Section \ref{s:trans} the perturbation variable
\be\label{pert2}
v(x,t):=u(x+\alpha(t),t)-\bar u(x)
\ee
for $u$ a solution of \eqref{cons}, where $\alpha$ is to be specified later
in a way appropriate for the task at hand.
Subtracting the equations for $u(x+\alpha(t), t)$ and $\bar u(x)$,
we obtain the nonlinear perturbation equation
\be\label{nlpert}
v_t-Lv= N(v)_x -\d_t \alpha (\phi+ \partial_x v),
\ee
where $L$ as in \eqref{L} denotes the linearized operator about $\bar u$, 
$\phi=\bar u_x$, and
\be\label{N}
N(v):=-(f(\bar u+ v)-f(\bar u)- df(\bar u)v)
\ee
where, so long as $|v|_{H^1}$ (hence $|v|_{L^\infty}$ and $|u|_{L^\infty}$) 
remains bounded, 
\ba\label{Nbds}
N(v)&=O(|v||v_x|),\\
\partial_x N(v)&=O(|v||v_{xx}|+|v_x|^2).\\
\ea

Recalling that $\d_x\bar u$ is a stationary
solution of the linearized equations $u_t=Lu$,
so that $L\d_x\bar =0$, or
$$
\int^\infty_{-\infty}G(x,t;y)\bar u_x(y)dy=e^{Lt}\bar u_x(x)
=\d_x\bar u(x),
$$
we have, applying Duhamel's principle to \eqref{nlpert},
$$
\begin{array}{l}
  \displaystyle{
  v(x,t)=\int^\infty_{-\infty}G(x,t;y)v_0(y)\,dy } \\
  \displaystyle{\qquad
  -\int^t_0 \int^\infty_{-\infty} G_y(x,t-s;y)
  (N(v)+\dot \alpha v ) (y,s)\,dy\,ds + \alpha (t)\d_x \bar u(x).}
\end{array}
$$
Defining 
\begin{equation}
 \begin{array}{l}
  \displaystyle{
  \alpha (t)=-\int^\infty_{-\infty}e(y,t) v_0(y)\,dy }\\
  \displaystyle{\qquad
  +\int^t_0\int^{+\infty}_{-\infty} e_{y}(y,t-s)(N(v)+
  \dot \alpha\, v)(y,s) dy ds, }
  \end{array}
 \label{alpha}
\end{equation}
following \cite{ZH,Z4,MaZ2,MaZ3}, 
where $e$ is defined as in \eqref{e}, 
and recalling the decomposition $G=E+ G_u+ \tilde G$ of \eqref{fulldecomp},
we obtain the {\it reduced equations}
\begin{equation}
\begin{array}{l}
 \displaystyle{
  v(x,t)=\int^\infty_{-\infty} (G_u+\tilde G)(x,t;y)v_0(y)\,dy }\\
 \displaystyle{\qquad
  -\int^t_0\int^\infty_{-\infty}(G_u+\tilde G)_y(x,t-s;y)(N(v)+
  \dot \alpha v)(y,s) dy \, ds, }
\end{array}
\label{v}
\end{equation}
and, differentiating (\ref{alpha}) with respect to $t$,
and observing that 
$e_y (y,s)\rightharpoondown 0$ as $s \to 0$, as the difference of 
approaching heat kernels,
\begin{equation}
 \begin{array}{l}
 \displaystyle{
  \dot \alpha (t)=-\int^\infty_{-\infty}e_t(y,t) v_0(y)\,dy }\\
 \displaystyle{\qquad
  +\int^t_0\int^{+\infty}_{-\infty} e_{yt}(y,t-s)(N(v)+
  \dot \alpha v)(y,s)\,dy\,ds. }
 \end{array}
\label{alphadot}
\end{equation}
\medskip

We emphasize that this (nonlocal in time) choice of $\alpha$ 
and the resulting reduced equations are different from those
of Section \ref{s:trans}.
As discussed further in \cite{Go,Z4,MaZ2,MaZ3,Z2},
$\alpha$ may be considered in the present context as defining a notion of
approximate shock location.

\subsection{Nonlinear damping estimate}

\begin{proposition}[\cite{MaZ3}]\label{damping}
Assuming (H0)-(H5), let $v_0\in H^{4}$, 
and suppose that for $0\le t\le T$, the $H^{4}$ norm of $v$
remains bounded by a sufficiently small constant, for $v$ as in
\eqref{pert2} and $u$ a solution of \eqref{cons}.
Then, for some constants $\theta_{1,2}>0$, for all $0\leq t\leq T$,
\begin{equation}\label{Ebounds}
|v(t)|_{H^4}^2 \leq C e^{-\theta_1 t} |v(0)|^2_{H^4} 
+ C \int_0^t e^{-\theta_2(t-s)} (|v|_{L^2}^2 + |\dot \alpha|^2) (s)\,ds.
\end{equation}
\end{proposition}

\begin{proof}
Energy estimates
essentially identical to those in the proof of Proposition \ref{damp1},
but involving the new forcing term $\dot{\alpha}(t) \d_x \bar u(x)$.
Observing that $\partial_x^j(\partial_x \bar{u})(x)=O(e^{-\eta|x|})$ 
is bounded in $L^1$ norm for $j\le 4$, we obtain (in the simpler,
uniformly elliptic case), taking account
of this new contribution, the inequality
\[
\partial_t |v|_{H^4}^2(t) \leq -\theta |\partial_x^{5} v|_{L^2}^2 +
C \left(|v|_{H^4}^2 + |\dot{\alpha}(t)|^2\right),
\]
$\theta>0$, for $C>0$ sufficiently large, so long as $|v|_{H^4}$ 
remains sufficiently small, yielding \eqref{Ebounds} as before
by Sobolev interpolation followed by Gronwall's inequality.
See also \cite{MaZ3,RZ}.
\end{proof}

\subsection{Short time existence theory}

\begin{lemma}[\cite{RZ}] \label{weightest}
Assuming (H0)-(H5), 
let $M_0:=|(1+|x|^2)^{3/4}v_0(x)|_{H^4}<\infty$, and suppose that, for $0\le t\le T$, the supremum of $ |\dot\alpha|$,
and the $H^{4}$ norm of $v$, determined by \eqref{nlpert}, each remain bounded by some constant $C>0$. Then there exists some $M=M(C)>0$ such that, for all $0\le t\le T$, 
\begin{equation}\label{weightedEbounds}
|(1+|x|^2)^{3/4}v(x,t)|_{H^4}^2 \le Me^{Mt} \big( M_0 +\int_0^t 
|\dot\alpha|^2(\tau)\, d\tau \big).
\end{equation}
\end{lemma}

\begin{proof}
This follows by standard Friedrichs symmetrizer estimates carried out in the weighted $H^4$ norm. Specifically, making the coordinate change 
$$
v=(1+|x|^2)^{3/4}w,
$$
we obtain from \eqref{nlpert} a similar equation
plus lower-order commutator terms,
and similarly in the equations for $\partial_x^j w$ for $j=1,\dots,4$. 
Performing the same energy estimates 
as carried out on \eqref{nlpert} in the proof of Lemma \ref{damping}, 
we readily obtain the result by Gronwall's inequality. 
We refer to \cite[Lemma~5.2]{RZ} for further details in the general 
partially parabolic case.  
\end{proof}

\begin{remark}\label{zetareg}
Using Sobolev embeddings and equation \eqref{nlpert}, we see that Lemma~\ref{weightest} immediately implies that, if $|(1+|x|^2)^{3/4}v_0(x)|_{H^4}<\infty$ and if $|v(\cdot,t)|_{H^4}$, $|\dot\alpha(t)|$ are uniformly bounded on $0\le t\le T$, then
\[
|(1+|x|^2)^{3/4}v(x,t)| \quad\mbox{and}\quad
|(1+|x|^2)^{3/4}v_t(x,t)|
\]
are uniformly bounded on $0\le t\le T$ as well.
\end{remark}

\subsection{Proof of nonlinear stability}

Decompose now the nonlinear perturbation $v$ as
\be\label{vdecomp}
v(x,t)=w(x,t)+z(x,t),
\ee
where
\be\label{wzdef}
w:=\Pi_{cs}v, \quad z:=\Pi_u v.
\ee
Applying $\Pi_{cs}$ to \eqref{v} and recalling commutator
relation \eqref{comm}, we obtain an equation
\ba \label{w}
  w(x,t)&=\int^\infty_{-\infty} \tilde G (x,t;y)w_0(y)\,dy \\
  &\quad -\int^t_0\int^\infty_{-\infty} \tilde G_y (x,t-s;y)
\tilde \Pi_{cs}(N(v)+
  \dot \alpha v)(y,s) dy \, ds
\ea
for the flow along the center stable manifold, parametrized by
$w\in \Sigma_{cs}$.

\bl\label{quadlem} Assuming (H0)--(H1), for $v$ lying initially
on the center stable manifold $\cM_{cs}$,
\be\label{zwbd}
|\d_x^r z(x,t)|\le Ce^{-\theta|x|}|w|_{H^2}^2
\ee 
for all $0\le r\le 4$, some $C>0$, 
so long as $|w|_{H^2}$ remains sufficiently small.
\el

\begin{proof}
%
By quadratic-order tangency at $\bar u$ of the center stable 
manifold to $\Sigma_{cs}$, estimate \eqref{tanbd1}, we
have immediately $|z|_{H^2}\le C|w|_{H^2}^2$, whence
\eqref{zwbd} follows by equivalence of norms for finite-dimensional
vector spaces, applied to the $p$-dimensional subspace $\Sigma_u$.
(Alternatively, we may see this by direct computation using
the explicit description of $\Pi_u v$ afforded by Lemma \ref{projlem}.)
\end{proof}

\begin{proof}[Proof of Theorem \ref{t:mainstab}]

Recalling by Theorem \ref{t:cstrans}
that solutions remaining for all time in a sufficiently
small radius neighborhood $\cN$ of the set of translates of $\bar u$
lie in the center stable manifold $\cM_{cs}$, we obtain trivially
that solutions not originating in $\cM_{cs}$ must exit $\cN$ in finite time,
verifying the final assertion of orbital instability with respect
to perturbations not in $\cM_{cs}$.

Consider now a solution $v \in \cM_{cs}$, or, equivalently,
a solution $w\in \Sigma_{cs}$ of \eqref{w} with 
$z=\Phi_{cs}(w)\in \Sigma_u$.
Define
\ba \label{zeta2}
 \zeta(t)&:= \sup_{y, 0\le s \le t}
 \Big( (|w|+|w_x|)(\theta+\psi_1+\psi_2)^{-1}(y,t)\\
&\quad 
+|w|_{H^4}(1+s)^{\frac{1}{4}} + |\dot \alpha (s)|(1+s) \Big).
\ea
We shall establish:

{\it Claim.} For all $t\ge 0$ for which a solution exists with
$\zeta$ uniformly bounded by some fixed, sufficiently small constant,
there holds
\begin{equation}
\label{claim}
\zeta(t) \leq C_2(E_0 + \zeta(t)^2).
\end{equation}
\medskip

{}From this result, provided $E_0 < 1/4C_2^2$, 
we have that $\zeta(t)\le 2C_2E_0$ implies
$\zeta(t)< 2C_2E_0$, and so we may conclude 
by continuous induction that
 \begin{equation}
 \label{bd}
  \zeta(t) < 2C_2E_0
 \end{equation}
for all $t\geq 0$, from which we readily obtain the stated bounds.
(By standard short-time $H^s$ existence theory together
with Remark \ref{zetareg}, $v\in H^4$ exists and $\zeta$ remains
continuous so long as $\zeta$ remains bounded by some uniform constant,
hence \eqref{bd} is an open condition.)
\medskip

{\it Proof of Claim.}
By \eqref{wtpbd}, $|(1+|x|^2)^{\frac{3}{4}}w_0|_{H^4}
= |(1+|x|^2)^{\frac{3}{4}}\Pi_{cs}v_0|_{H^4} \le CE_0$.
Likewise, by Lemma \ref{quadlem},
\eqref{zeta2}, \eqref{Nbds}, and \eqref{ptproj}, for $0\le s\le t$,
\ba\label{Nlast}
|\tilde \Pi_{cs}(N(v)+ \dot \alpha v)(y,s)|&\le C\zeta(t)^2 
\Psi(y,s)|,\\
|(N(v)+ \dot \alpha v)(y,s)|&\le C\zeta(t)^2 
\Psi(y,s)|,\\
\ea
$\Psi$ as in \eqref{source}, while
\be\label{dNlast}
|\Pi_{cs}\partial_y(N(v)+ \dot \alpha v)(y,s)|\le C\zeta(t)^2 
C(1+s)^{-\frac{1}{4}} (\theta+\psi_1+\psi_2)(x,t).
\ee

Combining \eqref{Nlast} and \eqref{dNlast} with representations
(\ref{w}) and (\ref{alphadot}) and
applying Lemmas \ref{iniconvolutions} and \ref{convolutions}, we obtain
 \ba\label{claimw}
  |w(x,t)| &\le
  \int^\infty_{-\infty} |\tilde G(x,t;y)| |w_0(y)|\,dy
   \\
 &\qquad +\int^t_0
  \int^\infty_{-\infty}|\tilde G_y(x,t-s;y)||\tilde \Pi_{cs}(N(v)+
  \dot \alpha v)(y,s)| dy \, ds \\
  & \le
  E_0 \int^\infty_{-\infty} |\tilde G(x,t;y)|(1+|y|)^{-3/2}\,dy
   \\
 &\quad +
C\zeta(t)^2 \int^t_0
  \int^\infty_{-\infty}|\tilde G_y(x,t-s;y)|
\Psi(y,s) dy \, ds \\
&\le
C(E_0+\zeta(t)^2)(\theta + \psi_1+\psi_2)(x,t)
\ea
and, similarly,
\ba\label{claimalpha}
 |\dot \alpha(t)| &\le \int^\infty_{-\infty}|e_t(y,t)|
  |u_0(y)|\,dy \\
  &\qquad +\int^t_0\int^{+\infty}_{-\infty} |e_{yt}(y,t-s)||\tilde \Pi_{cs}(N(v)+
 \dot \alpha v)(y,s)|\,dy\,ds\\
&\le \int^\infty_{-\infty}E_0 |e_t(y,t)|(1+|y|)^{-3/2}\, dy\\
&\qquad
+ \int^t_0\int^{+\infty}_{-\infty}C\zeta(t)^2 |e_{yt}(y,t-s)|\Psi(y,s) \,dy\,ds\\
&\le C(E_0+\zeta(t)^2)(1+t)^{-1}.\\
\ea

Likewise, we may estimate for $t\ge 1$
 \ba\label{claimwx}
  |w_x(x,t)| &\le
  \int^\infty_{-\infty} |\tilde G_x(x,t;y)| |w_0(y)|\,dy
   \\
 &\qquad +\int^{t-1}_0
  \int^\infty_{-\infty}|\tilde G_{xy}(x,t-s;y)||\tilde \Pi_{cs}(N(v)+
  \dot \alpha v)(y,s)| dy \, ds \\
 &\qquad +\int_{t-1}^t
  \int^\infty_{-\infty}|\tilde G_{x}(x,t-s;y)||\Pi_{cs}\partial_y (N(v)+
  \dot \alpha v)(y,s)| dy \, ds \\
  & \le
  E_0 \int^\infty_{-\infty} |\tilde G(x,t;y)|(1+|y|)^{-3/2}\,dy
   \\
 &\quad + C\zeta(t)^2 \int^{-1}_0 \int^\infty_{-\infty}|\tilde G_{xy}(x,t-s;y)| 
\Psi(y,s) dy \, ds \\
 &\quad + C\zeta(t)^2 \int^t_{t-1} \int^\infty_{-\infty}|\tilde G_x(x,t-s;y)| 
(\theta+\psi_1+\psi_2)(y,s) dy \, ds \\
&\le
C(E_0+\zeta(t)^2)(\theta + \psi_1+\psi_2)(x,t),
\ea
and for $0\le t\le 1$, substitute for the first term on the righthand
side instead
$$
  \int^\infty_{-\infty} |(\tilde G_x+\tilde G_y)(x,t;y)| |w_0(y)|\,dy
  +\int^\infty_{-\infty} |\tilde G(x,t;y)| |\partial_y w_0(y)|\,dy
$$
to obtain the same bound with the aid of \eqref{inicommeq}.

By Lemma \ref{quadlem}, 
\be\label{zh4}
|z|_{H^4}(t)\le C|w|_{H^2}^2(t)\le C\zeta(t)^2.
\ee
In particular, 
$
|z|_{L^2}(t)\le C\zeta(t)^2(1+t)^{-\frac{1}{2}}.
$
Applying Proposition \ref{damping} and using \eqref{claimw} and
\eqref{claimalpha}, we thus obtain
\be\label{claimwH4}
|w|_{H^4}(t)\le C(E_0+\zeta(t)^2)(1+t)^{-\frac{1}{4}}.
\ee
Combining \eqref{claimw}, \eqref{claimalpha}, \eqref{claimwx},
and \eqref{claimwH4}, 
we obtain \eqref{claim} as claimed.
\medbreak

As discussed earlier,
from \eqref{claim}, we obtain by continuous induction \eqref{bd}, or
$
 \zeta\le 2C_2|v_0|_{L^1\cap H^4}, 
$
whereupon the claimed bounds on $|v|$, $|v_x|$ and $|v|_{H^4}$ follow by
\eqref{claimw}, \eqref{claimwx} and \eqref{claimwH4} together with
\eqref{zwbd}, and the bounds on $|\dot \alpha|$ by \eqref{claimalpha}.
It remains to establish the bound on $|\alpha|$, expressing convergence
of phase $\alpha$ to a limiting value $\alpha(+\infty)$.

By Lemmas \ref{iniconvolutions}--\ref{convolutions}
together with the previously
obtained bounds \eqref{Nlast} and $\zeta\le CE_0$, 
and the definition \eqref{zeta2} of $\zeta$, the formal limit
$$
 \begin{aligned}
 \alpha(+\infty)&:=  \int^\infty_{-\infty}
e(y,+\infty) u_0(y)\,dy  \\
 &\qquad
 +\int^{+\infty}_0\int^{+\infty}_{-\infty} 
 e_y(y,+\infty)(N(v)+ \dot \alpha v)(y,s)\,dy\,ds \\
&\le 
 \int^\infty_{-\infty}
E_0 |e(y,+\infty)|(1+|y|)^{-3/2} \,dy  \\
 &\qquad
 +\int^{+\infty}_0\int^{+\infty}_{-\infty} 
CE_0 |e_y(y,+\infty)| \Psi(y,s) \,dy\,ds \\
&\le CE_0
\end{aligned}
$$
is well-defined, as the sum of absolutely convergent integrals.

Applying Lemmas \ref{iniconvolutions}--\ref{convolutions} a final time,
we obtain
$$
 \begin{aligned}
 |\alpha(t)-\alpha(+\infty)| &\le \int^\infty_{-\infty}
|e(y,t)-e(y,+\infty)| | v_0(y)| \,dy  \\
 &\qquad
 +\int^t_0\int^{+\infty}_{-\infty} |e_{y}(y,t-s)-e_y(y,+\infty)|\\
&\qquad \times
|(N(v)+ \dot \alpha v)(y,s)|\,dy\,ds \\
 &\qquad
 +\int_t^{+\infty}\int^{+\infty}_{-\infty} |e_y(y,+\infty)|
|(Q(v)+
  \dot \alpha v)(y,s)|\,dy\,ds \\
&\le \int^\infty_{-\infty}
E_0 |e(y,t)-e(y,+\infty)| (1+|y|)^{-3/2} \,dy  \\
 &\qquad
 +\int^t_0\int^{+\infty}_{-\infty} |e_{y}(y,t-s)-e_y(y,+\infty)|
CE_0 \Psi(y,s) \,dy\,ds \\
 &\qquad
 +\int_t^{+\infty}\int^{+\infty}_{-\infty} |e_y(y,+\infty)|
CE_0 \Psi(y,s)\,dy\,ds \\
&\le CE_0(1+t)^{-1/2},
\end{aligned}
$$
establishing the remaining bound and completing the proof.

\end{proof}

\medbreak
{\bf Acknowledgement.}
Thanks to Milena Stanislavova and Charles Li for two interesting discussions
that inspired this work, and to Milena Stanislavova
for pointing out the reference in \cite{GJLS}.



%
\appendix

\section{Proofs of miscellaneous lemmas}\label{extraproofs}

We include for completeness the proofs of earlier cited lemmas
that were not proved in the main body of the text.

\begin{proof}[Proof of Proposition \ref{linest}]
By sectoriality of $L$,
we have the inverse Laplace transform representations
\begin{equation}\label{ILT}
\begin{aligned}
 e^{tL}\Pi_u &:= \int_{{\Gamma_u }} e^{\lambda t} 
(\lambda - L)^{-1} \, d\lambda,\\
 e^{tL}\Pi_{cs} &:= \int_{{\Gamma_{cs}}} e^{\lambda t} 
(\lambda - L)^{-1} \, d\lambda,\\
\end{aligned}
\end{equation}
where $\Gamma_{cs}$ denotes a sectorial contour bounding the center and stable
spectrum to the right \cite{Pa}, which by (A2) may be taken
so that $\Re \Gamma_s\le \o$, and $\Gamma_u$ denotes a closed curve enclosing
the unstable spectrum of $L$, with $\Re \Gamma_s\ge \beta>0$.
Estimating $|e^{Lt}\Pi_j|_{L^2}\le \int_{\Gamma_j}|e^{\lambda t}|
|(\lambda-L)^{-1}|_{L^2} |d\lambda|$ using sectorial resolvent estimates,
we obtain the stated bounds from $L^2\to L^2$; see \cite{He}. 

Applying the resolvent formula 
$ L(\lambda-L)^{-1}= \lambda(\lambda-L)^{-1}- Id $,
we obtain 
$$
 Le^{tL}\Pi_j := \int_{{\Gamma_j}} e^{\lambda t} 
(\lambda - L)^{-1} \, d\lambda,
$$
from which we obtain $|Le^{Lt}\Pi_u|_{L^2}\le Ce^{-\beta t}$ 
for $t\le 0$ by boundedness of $\Gamma_u$, yielding the stated bound
from $L^2\to H^2$.
\end{proof}

\begin{proof}[Proof of Lemma \ref{projlem1}]
Recalling that $L$ has at most finitely many unstable eigenvalues,
we find that $\Pi_u$ may be expressed as
$$
\Pi_u f= \sum_{j=1}^p \phi_j(x) \langle \tilde \phi_j, f\rangle,
$$
where $\phi_j$, $j=1, \dots p$ are generalized right eigenfunctions of
$L$ associated with unstable eigenvalues $\lambda_j$, 
satisfying the generalized eigenvalue equation $(L-\lambda_j)^{r_j}\phi_j=0$,
$r_j\ge 1$, and $\tilde \phi_j$
are generalized left eigenfunctions.
Noting that $\phi_j$, $\tilde \phi_j$ and derivatives decay exponentially
by standard theory \cite{He,ZH,MaZ1}, and estimating 
$$
|\partial_x^j\Pi_u f|_{L^p}=|\sum_j \partial_x^j\phi_j \langle \tilde \phi_j f
\rangle|_{L^p}\le
\sum_j |\partial_x^j\phi_j|_{L^p} |\tilde \phi_j|_{L^q} |f|_{L^p}
\le C|f|_{L^p}
$$
for $1/p+1/q=1$, we obtain the claimed
bounds on $\Pi_u$, from which the bounds on
$\Pi_{cs}={\rm Id}-\Pi_u$ follow immediately.
\end{proof}

 \begin{proof}[Proof of Lemma \ref{trunc}]
 The norm in $H^2$ is a quadratic form, hence the map
 $$ v \in H^2 \mapsto \rho\Big( \frac{ | v |_{H^2}}{\delta}\Big) 
\in \R_+,$$
 is smooth, and $N^\delta$ is as regular as $N.$ 
The Lipschitz bound follows by
 \begin{eqnarray} | N^\delta(v_1) - 
N^\delta(v_2)|_{L^2} 
& \leq & | \rho\Big( \frac{ | v_1 |_{H^2}}{\delta}\Big) -   
 \rho\Big( \frac{ | v_2 |_{H^2}}{\delta}\Big) |_{L^\infty} 
| N(v_1)|_{L^2} \nonumber \\ &  & + \, | 
\rho\Big( \frac{ | v_2 |_{H^2}}{\delta}\Big) |_{L^\infty}  
| N(v_1) - N(v_2)|_{L^2} \nonumber \\
 & \leq & 3 | v_1 - v_2|_{H^2} 
\Big(
\sup_{|v|_{H^2} < \delta} \frac{| N(v)|_{L^2}}{\delta}
+   
\sup_{|v|_{H^2} < \delta} | dN(v)|_{L^2}
\Big) ,\nonumber
 \end{eqnarray}
 so that $\sup_{|v|_{H^2} < \delta}  | N(v)|_{L^2} = O( \delta^2),$
 $\sup_{|v|_{H^2} < \delta}  | dN(v)|_{L^2} = O( \delta).$
Finally, $|N^\delta(v)|_{L^2}\le |N(v)|_{L^2}\le C|v|_{H^2}^2$ 
for $|v|_{H^2}\le 2\delta$ 
by Moser's inequality, while $N^\delta(v)\equiv 0$ for $|v|_{H^2}\ge 2\delta$, 
yielding \eqref{Nquad}.
 \end{proof}

\begin{proof}[Proof of Proposition \ref{greenbounds}]
As observed in \cite{TZ2},
it is equivalent to establish decomposition 
\be\label{fulldecomp}
G=G_u + E+\tilde G
\ee
for the full Green function $G(x,t;y):=e^{Lt}\delta_y(x)$,
where 
$$
G_u(x,t;y):=\Pi_u e^{Lt}\delta_y(x)
=
e^{\gamma t}\sum_{j=1}^p\phi_j(x)\tilde \phi_j(y)^t
$$
for some constant matrix $M\in \C^{p\times p}$
denotes the Green kernel of the linearized solution operator
on $\Sigma_u$, $\phi_j$ and $\tilde\phi_j$ right and left
generalized eigenfunctions associated with unstable eigenvalues
$\lambda_j$, $j=1,\dots,p$.

The problem of describing the full Green function
has been treated in \cite{ZH, MaZ3}, 
starting with the Inverse Laplace Transform representation
\be\label{ILT2}
G(x,t;y)=e^{Lt}\delta_y(x)= \oint_\Gamma e^{\lambda t}(\lambda-L(\e))^{-1} 
\delta_y(x)d\lambda \, ,
\ee
where 
$$
\Gamma:= \partial \{ \lambda : \Re \lambda\le \eta_1 - \eta_2 |\Im \lambda|\}
$$
is an appropriate sectorial contour, $\eta_1$, $\eta_2>0$;
estimating the resolvent kernel 
$G^\eps_\lambda(x,y):=(\lambda-L(\e))^{-1}\delta_y(x)$
using Taylor expansion in $\lambda$,
asymptotic ODE techniques in $x$, $y$, and judicious decomposition
into various scattering, excited, and residual modes;
then, finally, estimating the contribution of various modes to \eqref{ILT2}
by Riemann saddlepoint (Stationary Phase) method, moving contour
$\Gamma$ to a optimal, ``minimax'' positions for each
mode, depending on the values of $(x,y,t)$.

In the present case, we may first move $\Gamma$ to a contour
$\Gamma'$ enclosing (to the left) all spectra of $L$
except for the $p$ unstable eigenvalues $\lambda_j$, $j=1, \dots, p$,
to obtain
$$
G(x,t;y)= \oint_{\Gamma'} e^{\lambda t}(\lambda-L)^{-1} d\lambda
+ \sum_{j=\pm} 
\Res_{\lambda_j(\eps)} \big( e^{\lambda t}(\lambda-L)^{-1}
\delta_y(x) \big),
$$
where
$\Res_{\lambda_j(\eps)} \big( e^{\lambda t}(\lambda-L)^{-1}
\delta_y(x) \big) = G_u(x,t;y)$, then estimate the remaining term
$ \oint_{\Gamma'} e^{\lambda t}(\lambda-L)^{-1} d\lambda$
on minimax contours as just described.
See the proof of Proposition 7.1, \cite{MaZ3}, for a detailed 
discussion of minimax estimates $E+G$ and of Proposition 7.7, 
\cite{MaZ3} for a complementary discussion of residues
incurred at eigenvalues in $\{\Re \lambda\ge 0\}\setminus\{0\}$.
See also \cite{TZ2}.
\end{proof}

\end{document}